\documentclass[hidelinks,onefignum,onetabnum]{siamart251216}

\usepackage{graphicx} 
\usepackage{amsmath,amssymb,bm,url}
\usepackage{hyperref,cleveref}
\usepackage{graphicx}
\usepackage{caption}
\usepackage{subcaption}
\usepackage[T1]{fontenc}
\usepackage{comment}
\usepackage{booktabs}
\usepackage{xcolor}
\usepackage{array}

\usepackage{cleveref}
\usepackage{algorithm}
\usepackage{algorithmicx}
\usepackage{algpseudocode}
\usepackage{hyperref}
\usepackage{cancel}

\usepackage{booktabs}
\usepackage{longtable}
\usepackage{geometry}
 \usepackage{enumitem}

\algrenewcommand\algorithmicrequire{\textbf{Input:}}
\algrenewcommand\algorithmicensure{\textbf{Output:}}

\newcommand{\vv}{\bm{v}}
\newcommand{\vx}{\textbf{x}}

\newcommand{\R}{\mathbb{R}}
\newcommand{\di}{\,\mathrm{d}}

\newcommand{\x}{\mathbf{x}}
\newcommand{\z}{\mathbf{z}}

\newcommand{\M}{\mathcal{M}}

\newcommand{\E}[1]{\mathbb{E}\left[#1\right]}

\newcommand{\br}[1]{\left(#1\right)}

\newcommand{\abs}[1]{\left\vert#1\right\vert}

\newcommand{\bbr}[1]{\left\{#1\right\}}

\newcommand{\vbr}[2]{#1(#2 #1)}
\newcommand{\rbr}[1]{\left[#1\right]}

\newcommand{\B}[1]{\boldsymbol{#1}}

\newcommand{\cb}[1]{{\color{blue}#1}}

\makeatletter
\renewcommand\tableofcontents{%
  \@starttoc{toc}%
}
\newcommand*\l@section[2]{%
  \ifnum \c@tocdepth >\z@
  \addpenalty\@secpenalty
  \addvspace{1.0em \@plus\p@}%
  \setlength\@tempdima{1.5em}%
  \begingroup
    \parindent \z@ \rightskip \@pnumwidth
    \parfillskip -\@pnumwidth
    \leavevmode \bfseries
    \advance\leftskip\@tempdima
    \hskip -\leftskip
    #1\nobreak\hfil \nobreak\hb@xt@\@pnumwidth{\hss #2}\par
  \endgroup
  \fi}
\newcommand*\l@subsection{\@dottedtocline{2}{1.5em}{2.3em}}
\newcommand*\l@subsubsection{\@dottedtocline{3}{3.8em}{3.2em}}
\makeatother

\begin{document}

\title{Deep Neural networks for solving 
high-dimensional parabolic partial differential equations}
\author{Wenzhong Zhang \thanks{Suzhou Institute for Advanced Research, University of Science and Technology of China, Suzhou, Jiangsu, 215127, P. R. China (current address), and  Department of Mathematics, Southern
Methodist University, Dallas, TX 75275, USA (primary address).}
\and Zheyuan Hu \thanks{Joint first author. Department of Computer Science, School of Computing, National University of Singapore, 119077, Singapore.}
\and Wei Cai\thanks{Corresponding author (\texttt{cai@smu.edu}), Department of Mathematics, Southern
Methodist University (SMU), Dallas, TX 75275, USA. This author is supported by the Clements Chair in Applied Math at SMU.}
\and George EM Karniadakis\thanks{Division of Applied Mathematics, Brown University, Providence, RI 02912, USA. This author would like to acknowledge  the MURI/AFOSR FA9550-20-1-0358 project, and the DOE-MMICS SEA-CROGS DE-SC0023191 awards}.}

\maketitle

\begin{abstract}
The numerical solution of high-dimensional partial differential equations (PDEs) is severely constrained by the curse of dimensionality (CoD), rendering classical grid-based methods impractical beyond a few dimensions. In recent years, deep neural networks have emerged as a promising mesh-free alternative, enabling the approximation of PDE solutions in tens to thousands of dimensions. This review provides a tutorial-oriented introduction to neural-network–based methods for solving high-dimensional parabolic PDEs, emphasizing conceptual clarity and methodological connections.
We organize the  literature around three unifying paradigms: (i) PDE-residual–based approaches, including physics-informed neural networks and their high-dimensional variants; (ii) stochastic methods derived from Feynman–Kac and backward stochastic differential equation formulations; and (iii) hybrid derivative-free random difference approaches designed to alleviate the computational cost of  derivatives in high dimensions. For each paradigm, we outline the underlying mathematical formulation, algorithmic implementation, and practical strengths and limitations.
Representative benchmark problems -- including Hamilton-Jacobi-Bellman and Black-Scholes equations in up to 1000 dimensions -- illustrate the scalability, effectiveness, and accuracy of the methods. The paper concludes with a discussion of open challenges and future directions for reliable and scalable solvers of high-dimensional PDEs.

\end{abstract}



\setcounter{tocdepth}{2}
\tableofcontents
 
\section{Introduction}

The numerical solution of high-dimensional partial differential equations (PDEs) remains a fundamental challenge in scientific computing due to the 
curse-of-dimensionality (CoD), a term introduced by Bellman in the context of dynamic programming \cite{bellman1965dynamic}. As the number of dimensions 
$d$ increases, the computational complexity of conventional discretization methods such as finite differences or finite elements typically scales as 
$O(N^d)$, where 
$N$ is the number of degrees of freedom per dimension. This exponential growth renders classical grid-based methods infeasible for even moderately large $d$. Over the past two decades, a variety of sparse approximation and low-rank methods have been proposed to mitigate CoD, including sparse grids \cite{bungartz2004sparse}, low-rank tensor methods \cite{khoromskij2014tensor}, and the use of polynomial chaos and stochastic Galerkin expansions  with sparse grids for uncertainty quantification \cite{xiu2002wiener}, and multi-level Monte Carlo Picard iterations \cite{Jentzen21}.  More recently, deep neural networks have emerged as promising function approximators for high-dimensional PDEs by exploiting their ability to capture complex solution structures. 
Approaches such as DeepBSDE \cite{han2016deep,Han2018DeepBSDE},  DGM \cite{sirignano2018dgm}, Physics-Informed Neural Networks (PINNs) \cite{raissi2019pinn}  in high-dimensional variants \cite{hu2023tackling,hu2023hutchinson,he2023learning}, Deep Ritz methods \cite{Weinan2017TheDR}, DeepMartnet \cite{Cai2025SOCMart,Cai2025Mart}, and Deep Picard Iteration \cite{han2024picard} aim to address the CoD by reformulating the PDE solution as an optimization problem in a high-capacity neural network function space. These methods offer scalable alternatives for problems in which traditional methods become intractable, particularly in application domains such as kinetic theory, financial mathematics, stochastic optimal controls, uncertainty quantification, biological systems, and quantum physics.

\label{figure:cod}

High-dimensional PDEs are ubiquitous in science and engineering, solutions to which can lead to numerous practical impacts and contributions.
The Schrödinger equation is the PDE in quantum physics and chemistry that describes molecular interactions and applies to molecular property prediction and new drug/material discovery. The PDE dimension equals the number of atoms in a molecule or a nano-system, which can be high in large molecules and material samples.
High-dimensional PDEs also apply to economics/econometrics, namely, the monetary and fiscal policies can be formulated as a stochastic optimal control problem, satisfying the Hamilton-Jacobi-Bellman (HJB) equation. More specifically, the HJB equation determines how to control these policies to minimize unemployment,  management cost, and maximize growth and profit optimally, where the PDE dimension equals the number of policies the government can take.
The Black-Scholes (BS) equation for option pricing in financial engineering is also high-dimensional, whose dimension equals the number of stocks to model the complex correlation between them.
Moreover, in the fast emerging area of data sciences,  state-of-the-art generative models like diffusion models \cite{song2021scorebased} are described by stochastic differential equations (SDEs), where the probability density function (PDF) of a generative model satisfies a Fokker-Planck (FP) equation. The dimension of this PDE equals the image dimension, whose scale is $\sim10^4-10^5$ in high-resolution image synthesis.

In this paper, we will present neural network methods designed to mitigate the difficulties posed by the CoD in solving, in particular, the following  $d$-dimensional quasi-linear parabolic PDE,
\begin{equation}
\partial_{t} u+\mathcal{L}(u)=f(t,\textbf{x}, u), \quad \textbf{x} \in \Omega  \label{eq-paraPDE1}%
\end{equation}
where $\Omega$ could be a bounded domain or the whole space  $\mathbb{R}^d$. For the bounded domain case,  a boundary condition (BC) will be imposed
\begin{equation}
    \mathcal{B}u(t,\vx)=b(t,\vx) \quad \text{on}\ \Gamma =\partial \Omega, \quad t>0,
    \label{BC}
\end{equation}
and, a decaying condition at infinity for the case of the whole space. In addition, the following terminal condition at $t=T$ will be included,
\begin{equation}
u(T,\textbf{x})=g(\textbf{x}).
\label{termData}
\end{equation}
Here, the quasi-linear differential operator  $\mathcal{L}$ is  defined by
\begin{equation}
\mathcal{L}=\mu(t,\textbf{x}, u) \cdot \nabla + \frac{1}{2} \sum_{i,j=1}^d a_{i,j}(t,\textbf{x}, u) \frac{\partial^2}{\partial x_i \partial x_j} =\mu(t,\textbf{x}, u) \cdot \nabla +\frac{1}{2}\mathrm{Tr}(A\nabla \nabla^\intercal),
\label{generator_ndim}
\end{equation}
where the $d \times d$ is a positive definite  coefficient matrix $A(t.\textbf{x}, u)=(a_{i,j}(t,\textbf{x}, u))$. Here, we consider the backward parabolic PDE with a terminal condition as often encountered in control problems.
By setting a new time $s=T-t$, the terminal value problem of the parabolic PDE becomes the conventional initial value problem for the time variable $s$ for a parabolic PDE such as the Fokker-Planck equation. 
 
For the boundary value problem of elliptic PDEs, the time derivative operator $\partial_t$ will be dropped in \eqref{eq-paraPDE1} as well as the terminal condition \eqref{termData} for a time-independent solution $u(\vx)$, while the boundary condition becomes
 \begin{equation}
    \mathcal{B}u(\vx)=b(\vx) \quad \text{on}\ \Gamma =\partial \Omega.
    \label{ellipticBC}
\end{equation}


Due to the compositional structure with affine mapping and nonlinear activation function  in the deep neural networks (DNNs), they can be used to represent high-dimensional functions efficiently,  much research is though needed to understand why the DNN could provide good approximations to high dimensional functions, which will not be addressed in this paper (see \cite{WeinanToward20,WojWeinan21}).
Meanwhile, the development of using the DNNs for solving high-dimensional parabolic PDEs has shown their promising potential in addressing the CoD problems, which are beyond the reach of traditional numerical methods. It is the objective of this paper to introduce some recent developments in DNNs so that researchers, especially graduate students, new to this field can get a good knowledge of some of the approaches being used in mitigating the effect of CoD in solving high dimensional PDEs and their applications in various fields in sciences, engineering, and economics. A more comprehensive overview in this field can be found in  \cite{HanEreview21}.

\section {Conceptual strategies for overcoming CoD}
The neural network methods to be discussed in this tutorial employ three main strategies:
\smallskip
\begin{itemize}[leftmargin=1em]
    \item {\bf Strategy A: PDE residual based PINNs - Monte Carlo in physical space}. In this approach, the PDE based NNs uses the residual of the PDEs for the NNs as the loss function to satisfy the PDEs at the Monte Carlo (MC) sampled locations in the solution domain, and the training aims to drive the NN's PDE residual to zero to ensure the neural network to satisfy the PDE.
    \item {\bf Strategy B: SDE-based NN - Monte Carlo in path space}. Stochastic differential equation (SDE) based NNs build the loss function using the Pardoux-Peng nonlinear Feynman-Kac formula for quasi-linear parabolic and elliptic PDEs, which avoids computing the Hessian matrix, and, the training objective is to enforce the validity of the probabilistic representation of the PDE solution using SDE paths.
    \item {\bf Strategy C: Random differences by expectation - operator estimation without differentiations}. In this approach, which embodies some aspects of the two strategies above, random difference operators based on SDE representation of the PDEs solution are derived for the PDE differential operator avoiding using any derivative operations while enforcing PDEs residual conditions by a strong form (PINN) or a Galerkin weak form (under a weight given by the path probability density function in the path space).
\end{itemize}

\subsection{Strategy A: PDE residual based PINN methods}

We consider the following elliptic PDE version of \eqref{eq-paraPDE1} (without the time derivative) defined on a domain $\Omega \subset \mathbb{R}^d$ with the boundary condition \eqref{ellipticBC} defined on $\Gamma \subset \mathbb{R}^d$,
\begin{equation}\label{eq:PDE}
\begin{aligned}
\mathcal{L}u(\vx)=f(\vx) \ \text{in}\ \Omega. 
\end{aligned}
\end{equation}

PINN \cite{raissi2019pinn} is a neural network-based PDE solver via minimizing the following boundary $\mathcal{L}_b(\theta)$ and residual loss $\mathcal{L}_r(\theta)$ functions,
\begin{equation}
\begin{aligned}
\mathcal{L}(\theta) &= \lambda_b \mathcal{L}_b(\theta) + \lambda_r \mathcal{L}_r(\theta)\\
&=\frac{\lambda_b}{n_b}\sum_{i=1}^{n_b} {|\mathcal{B}u_{\theta}(\vx_{b,i})-b(\vx_{b,i})|}^2 + \frac{\lambda_r}{n_r}\sum_{i=1}^{n_r} {|\mathcal{L}u_{\theta}(\vx_{r,i})-f(\vx_{r,i})|}^2, 
\end{aligned}
\label{pinnloss}
\end{equation}
where $\lambda_b, \lambda_r > 0$ are the weights used to balance the losses. 
The PINN loss is computed on randomly sampled points in the physical space as in traditional Monte Carlo methods but does not rely on a grid and mesh as in traditional grid based methods. Hence, PINN is promising in tackling
the CoD. 


%

\subsubsection{Methodology, architecture, and algorithm in PINNs}
Wang et al. \cite{wang20222} discovered that learning Hamilton-Jacobi-Bellman (HJB) equations in high dimensions in optimal control requires the use of the $L^\infty$
loss. 
Separable PINN \cite{cho2022separable} adopts a per-axis sampling strategy for residual points in high-dimensional spaces, which reduces computational cost and enables larger batch sizes through tensor product construction. However, this approach is primarily designed to accelerate training in low-dimensional settings (typically fewer than four dimensions) and targets the increase of collocation points in 3D PDEs, with its sampling point separation becoming impractical in truly high-dimensional problems.
Wang et al. \cite{wang2022tensor, wang2022solving} introduced tensor neural networks incorporating efficient numerical integration and separable architectures to address high-dimensional Schrödinger equations arising in quantum physics.
Zhang et al. \cite{zhang2020learning} tackled stochastic differential equations by employing PINNs in conjunction with spectral dynamically orthogonal and bi-orthogonal expansions to represent the solution space. 

As PINNs can be viewed as following traditional collocation methods for solving PDEs, DNNs mimicking classic finite element methods such as Galerkin and Ritz finite element methods have also been developed.
Zang et al. \cite{Zang2020Weak} developed a weak adversarial network framework that treats PDEs through their variational (weak) form.
The deep Galerkin method (DGM) \cite{sirignano2018dgm} addresses high-dimensional PDEs by training NNs to fulfill the governing equations by using PDE residual and Monte Carlo techniques to approximate derivatives. The deep Ritz method (DRM) \cite{Weinan2017TheDR} formulates the solution of high-dimensional PDEs as an energy minimization problem optimized through neural networks.

\subsubsection{Mitigation of auto-differentiation (AD) cost in PINNs}
As computing the residual loss involving high-order and high-dimensional neural networks' derivatives is costly, multiple amortizations have been proposed. 
Stochastic dimension gradient descent (SDGD) \cite{hu2023tackling} mitigates the cost via decomposing the high-dimensional PDE operator along each dimension. Part of all dimensions is sampled and forms an unbiased estimator of the full operator. 
Hutchinson trace estimation (HTE) \cite{hu2023hutchinson} treats the ubiquitous second-order PDE operator as the Hessian trace and proposes using it for randomized, unbiased estimation. 
%
%
Stochastic Taylor derivative estimator (STDE) \cite{shi2024stde} extends SDGD and HTE to general PDEs and provides a comprehensive illustration and comparison of forward,
backward, and Taylor modes AD in PyTorch \cite{paszke2019pytorch} or JAX \cite{jax2018github}.
He et al. \cite{he2023learning,hu2023bias} proposed parameterizing PINNs using a Gaussian smoothed model and optimizing them via Stein's identity, thereby eliminating the need for back-propagation and reducing the extensive differentiation typically involved.
%
%
Additional efforts have focused on designing special forward computation rules tailored to specific Laplacian operators \cite{li2024computational,li2024dof}.

\subsection{Strategy B: SDE-based NNs}

Solution to the quasi-linear parabolic PDE \cref{eq-paraPDE1} 
has a probabilistic representation by the Feynman--Kac formula, or by the more general nonlinear Feynman--Kac formula in the Pardoux--Peng theory of forward-backward SDE (FBSDE) formulation \cite{Pardoux1990FBSDE} (refer to Section 4.1 for more details).
The FBSDE consists of a forward SDE of stochastic process $X_t$, a backward SDE of stochastic process $Y_t$, as well as an auxiliary process $Z_t$.
Solution to the FBSDEs is connected to that of the PDE \cref{eq-paraPDE1} by
\begin{equation*}
    Y_t = u(t, X_t), \quad Z_t = \nabla u(t, X_t),
\end{equation*}
see \cref{SDE-PDE}.
Major benefits of the FBSDE formulation include the fact that exploring the domain along trajectories of the stochastic process $X_t$ requires only time discretization, thus preventing the CoD of data, and that evaluating second-order derivatives in the differential operator $\mathcal{L}$ of the PDE is avoided.



We roughly categorize the FBSDE based methods for solving quasi-linear parabolic PDEs \cref{eq-paraPDE1} by the way the FBSDE formulation \cref{bsde1} is utilized:  (1) Backward schemes that recursively reduce to sub-problems on a smaller time domain; (2)  SDE path loss schemes where loss functions evaluated at the terminal time or along trajectories of $X_t$ play an important role.


\subsubsection{Backward schemes}
Solution to the process $Y_t$ provided with terminal values can be progressively solved in a dynamical programming (DP) manner, i.e. in backward order of time steps, so that the problem is also progressively resolved.
When stepwise DNNs are trained, transfer learning can be applied from the outcome of previous steps to newly initialized networks.
The Feynman--Kac formula \cite{weinanbook} represents $Y_t$ by an expectation conditioned by the natural filtration at time $t$.
For quasi-linear PDEs, the evaluation often involves future values of $Y_s$ for $s \in (t, T]$, which may be handled stepwisely to reduce approximation error.
Huré et al. \cite{Pham2020DBDP} proposed backward DP-type schemes training DNNs that approximate $u(t, \cdot)$ and $\nabla u(t, \cdot)$ on descending time steps, with the loss functions defined as square of the residual of discretized SDE on each time step.
Pham et al. \cite{pham2021neural} extended their work  to fully nonlinear PDEs.
Beck et al. \cite{Beck2021DeepSplitting} developed backward schemes for $u(t, \cdot)$ that are ruled by the $L^2$-minimality property of stepwise conditional expectations \cite{klenke2008probability}. As the Feynman-Kac formula can be viewed as a generalized characteristic method connecting the terminal condition to the solution at earlier times (backward).  Han et al. \cite{han2024picard}  incorporates Monte Carlo simulations of the Feynman--Kac formula \cite{weinanbook} for $Y_t$  and the Malliavin method \cite{ma2002representation} for $Z_t$ for a Picard fix-point iteration scheme for the PDE solution represented by a DNN.

\subsubsection{SDE path loss schemes}
The Deep BSDE method \cite{han2016deep,Han2018DeepBSDE} by Han et al. is the earliest known deep learning method for  solving high dimensional PDEs.
The method finds $Y_0$ from FBSDE to provide approximation of PDE solution at $t=0$ by matching simulated $Y_T$ with the terminal condition of PDE, where DNNs approximating the auxiliary process $\sigma(t, X_t)^{\top} Z_t$ on each time step of SDE simulation are employed.
Also, two different representations of $Y_t$ from the FBSDE and $u(t, X_t)$ from the PDE for any $t \in [0, T]$ inspire deep learning methods that exploit information along the trajectories of $X_t$ and $Y_t$.
Raissi proposed the FBSNN method \cite{Raissi2023FBSNNs} that uses a single DNN to approximate the PDE solution $u(\cdot, \cdot)$, and uses the difference between $Y_t$ from the discretized SDE and $u(t, X_t)$ from the DNN along trajectories of $Y_t$ as part of the loss function to provide approximation to solution along paths.
Based on this idea, Zhang and Cai \cite{Zhang2022FBNN} developed modified schemes that match solutions to different SDEs in the continuous-time limit, giving an order of $\mathcal{O}(\Delta t^{1/2})$ accuracy in the solution.
Kapllani and Teng \cite{kapllani2025backward} applied the Malliavin method to improve gradient approximations of FBSNNs.
Nüsken and Richter \cite{Nusken2023JML} proposed a diffusion loss scheme that utilizes the stochastic integral form of the backward SDE to build up the loss function, where the length of time interval for integration lies between a global limit as the Deep BSDE method, and a local limit as the PINN method. In \cite{HanNica2020}, the Feynman-Kac formula was used to build a derivative-free loss function for a PINN-type method. For HJB equations which can be interpreted as a stochastic optimal control (SOC) problem with PDE solution's derivative serving as the control, Zhou et al. \cite{Zhou2021Actor} developed temporal difference and policy improvement methods in the actor-critic framework,  where loss functions for the value function employ both the BSDE and Feynman-Kac formula for the value function.
Cai et al. \cite{Cai2023DeepMart,Cai2025SOCMart} proposed a martingale based approach, DeepMartNet, where the loss function enforces the martingale condition of the PDEs solution in  the SDE path spaces, in a Galerkin-type  adversarial learning framework. In \cite{DFMC20}, a diffusion Monte Carlo method using the FBSDEs was proposed for eigenvalue problems.

\subsection{Strategy C: Random differences by expectation}
When simulation of $Y_t$ is implemented with DNNs, the DNN calls accumulate under (usually nonlinear) operations in the scheme, slowing down the backward propagation, thus limiting the scalability for solving very high-dimensional problems or on a large time interval.
Such issue may be mitigated by slicing the trajectory into short pieces if localized properties of trajectories can be effectively used.
Xu and Zhang \cite{Xu2025Shotgun} showed that a single-step SDE discretization by FBSNN, giving a random difference operator, is an estimator of the residual of the PDE, and proposed a PINN-based scheme that applies coarse time sampling for accessing the approximated distribution of $X_t$.
Cai et al. \cite{Cai2025RDO} proposed a Galerkin weak form  of the random difference approach, an equivalent to DeepMartNet.

\section{PINNs for high-dimensional PDEs}


{\color{black}{\bf Sampling strategy}. The residual points in PINNs are often selected either along the trajectories of the stochastic processes related to the PDE operator or from a simple distribution that approximates those trajectories. Traditional residual-adaptive sampling or uniform sampling within a bounded domain unfortunately does not work if the PDE is defined over a high-dimensional, unbounded domain.}

In high dimensions, PINN's major computational cost comes from the residual loss involving evaluating the linear differential operator in equation (\ref{generator_ndim}).
This is a nontrivial task. For  high-dimensional PDEs, the PINN network size grows, and the amount of AD required for computing residual loss increases, leading to a magnified cost. Given a specific computational budget, there exists a dimension high enough such that even using one collocation point can lead to insufficient memory. This  cost cannot be further reduced since the batch size over the residual collection point is already 1, and it is unsolvable by traditional parallel computing since the minimal batch size should be 1.
Several mitigation procedures have been proposed to avoid such costly calculations.

\textbf{SDGD} \cite{hu2023tackling} decomposes a linear differential operator along each dimension 
\vspace*{-0.3cm}
\begin{equation}
\mathcal{L}= \sum_{i=1}^d \left(\sum_{j=1}^d a_{i,j}(t,\textbf{x}) \frac{\partial^2}{\partial x_i \partial x_j}\right) = \sum_{i=1}^d \mathcal{L}_{i}, \quad \mathcal{L}_{i}:= \sum_{j=1}^d a_{i,j}(t,\textbf{x}) \frac{\partial^2}{\partial x_i \partial x_j}.
\end{equation}
A random dimension index set $I \subset \{1,2,\cdots,d\}$ is sampled and the unbiased estimated residual of a PINN network $u_\theta$ is
\begin{equation}
\mathcal{L}u_\theta = \sum_{i=1}^d \mathcal{L}_{i}u_\theta \approx \frac{d}{|I|}\sum_{i \in I}\mathcal{L}_{i}u_\theta,
\end{equation}
where we reduce the entire Hessian to dimension pieces or several Hessian rows.

Instead of optimizing the expensive full PINN residual loss
\begin{equation}
\mathcal{L}_{\text{PINN}}(\theta) = \frac{1}{2}\left(\sum_{i=1}^d \mathcal{L}_iu_\theta(\vx) - f(\vx)\right)^2,
\end{equation}
we plug the estimated residual into the residual loss with a randomly sampled index set $I \subset \{1,2,\cdots,d\}$:
\begin{equation}
\mathcal{L}_{\text{SDGD}}(\theta) = \frac{1}{2}\left(\textcolor{black}{\frac{d}{|I|}\sum_{i \in I}} \mathcal{L}_iu_\theta(\vx) - f(\vx)\right)^2. 
\end{equation}
Due to the nonlinearity of the mean square error loss, it is biased. We can sample the index set twice with $I, J \subset \{1,2,\cdots,d\}$ independently to debias, as detailed in Appendix  A.3.

In sum, SDGD speeds up and scales up PINN in high dimensions via random dimension samples, which unbiasedly estimates the exact residual.
By minimizing the SDGD's estimated stochastic residual, we minimize the true residual loss in expectation.
SDGD can be combined with modern optimizers like Adam and SOAP.
SDGD's philosophy is that the algorithm's gradient can be decomposed based on different dimensions. Then, a subset of dimensional pieces can be sampled to construct an unbiased stochastic gradient for optimization. This avoids the computationally expensive calculation of the full gradient.
We sample a proper number of residual points in the low-dimensional PINN rather than an excessive number. This is to reduce the cost of each gradient while keeping it unbiased for convergence. When it comes to high dimensions, SDGD samples dimensions just like we sample points in traditional stochastic gradient descent (SGD).

\textbf{HTE} \cite{hu2023hutchinson} views the PDE operator as Hessian trace and uses Hutchinson trace estimation \cite{hutchinson1989stochastic} to mitigate the cost:
\begin{align}
\operatorname{Tr}(A) = \mathbb{E}_{\vv \sim p(\vv)}\left[ \vv^\intercal A \vv\right],
\end{align}
for all random variable $\vv \in \mathbb{R}^d$ such that $\mathbb{E}_{\vv \sim p(\vv)} [\vv\vv^T] = I$.
Therefore, the trace can be estimated by a Monte Carlo method:
\begin{equation}
\operatorname{Tr}(A)\approx \frac{1}{V}\sum_{i=1}^V \vv_i^\intercal A\vv_i,
\end{equation}
where each $\vv_i\in\mathbb{R}^d$ are $i.i.d.$ samples from $p(\vv)$.
Classically, to minimize the variance of HTE, we opt for the Rademacher distribution as follows: for each dimension of the vector $\vv \sim p(\vv)$, it is a discrete probability distribution that has a 50\% chance of getting +1 and a 50\% chance of getting -1.

Similar to SDGD, HTE replaces the original expensive PINN residual with a randomized HTE estimator.
Then, we plug the HTE estimator into the PINN loss as detailed in Appendix A.3.
%
HTE reduces the full Hessian to a Hessian-Vector Product (HVP).
The full Hessian is a huge $d \times d$ matrix, while HVP is just a scalar, reducing memory cost significantly.
The HVP of a PINN can be computed efficiently via Taylor mode \cite{shi2024stde,bettencourt2019taylormode} auto-differentiation in JAX.

\textbf{Connection between SDGD and HTE}. SDGD can be viewed as HTE with a specific choice of $p(\vv)$.
Consider the discrete distribution $p(\vv)$ such that $\vv = \sqrt{d}\boldsymbol{e}_i$ for $i=1,2,\cdots,d$ with probability $1 / d$, where $\boldsymbol{e}_i$ denotes the standard basis vector with a one in the $i$th coordinate and zero elsewhere.
Then, $\mathbb{E}_{\vv \sim p(\vv)}[\vv\vv^\mathrm{T}] = I$, and the corresponding estimator is $\operatorname{Tr}(A) \approx \frac{d}{|I|}\sum_{i\in I}A_{ii}$ where $I \subset \{1,2,\cdots,d\}$ is a multiset and $|I|$ denotes its cardinality, which is exactly sampling dimensionality as SDGD.
The only difference is that SDGD samples dimension without replacement. Hence, $I$ is a set rather than a multiset.
This view enables comparing their variance: SDGD from diagonal while HTE from off-diagonal, as detailed in Appendix A.3.

\textbf{STDE} \cite{shi2024stde} extends SDGD and HTE to general and arbitrary PDEs. It also provides a comprehensive illustration and comparison of forward, backward, and Taylor modes AD, as well as JAX advantage over PyTorch in AD.

\textbf{Randomized Smoothing (RS) PINN} \cite{he2023learning,hu2023bias} considers the following RS neural network structure for backpropagation-free PINN's AD computation:
$
u_\theta(\vx) = \mathbb{E}_{\delta \sim \mathcal{N}(0, \sigma^2\boldsymbol{I})}f_\theta(\vx+\delta),
$
where $f_\theta(\vx+\delta)$ is a vanilla neural network parameterized by $\theta$; $u_\theta(\vx)$ is the corresponding RS-PINN model serving as the surrogate model in fitting PDE solutions.
Given the special parameterization, $u_\theta(\vx)$'s derivatives can be analytically computed without backpropagation as detailed in Appendix A.3.

\section{SDE-based neural networks for quasi-linear PDEs} 
Two main advantages of using SDEs are: (1) avoiding the calculation of high-dimensional Hessian matrix of the PDE solution; and (2) SDE paths provide a natural sampling strategy of the solution domain.
Throughout this section, we consider the probabily space $(\Omega, \mathcal{F}, \mathbb{P})$ equipped with a natrural filtration $\mathbb{F} := \{\mathcal{F}_t\}_{0\le t\le T}$ generated by a d-dimensional standard Brownian motion $B_t$.

\subsection{It\^o formula and Pardoux-Peng nonlinear Feynman-Kac formula}\label{sect_pardoux_peng}

Let $X:[0,T]\times\Omega\to\R^d$ be an It\^{o} diffusion process; that is, $X$ is an $\mathbb{F}$-adapted process satisfying the forward stochastic differential equation
\begin{equation}
dX_t = \mu(t, X_t)dt+\sigma(t, X_t)dB_t,
\end{equation}
where the initial state $X_0$ is given, and
the functions $\mu: [0, T] \times \R^d \to \R^d$ and  $\sigma: [0, T] \times \R^d \to \R^{d \times d}$ are the drift and diffusion coefficients, respectively.
Let $h: [0, T] \times \R^d\to\R$ be a function of class $C^{1,2}$, i.e., continuously differentiable in $t$ and twice continuously differentiable in $\x$.
Then It\^{o}'s formula \cite{weinanbook} says that the process $h(t, X_{t})$ is also an It\^{o} process, governed by the forward SDE
\begin{equation}
dh(t, X_{t})=\left(\frac{\partial}{\partial t}+\mathcal{L}\right)h(t,X_{t})dt+\nabla^{\top} h(t,X_{t})\sigma
(t,X_t)dB_t, \label{ItoForm2}%
\end{equation}
in terms of the generator $\mathcal{L}$   for the diffusion $X_{t}$  given in \eqref{generator_ndim}
where the $d \times d$ diffusion coefficient matrix is 
\begin{equation}
    A(t,\textbf{x})=(a_{i,j}(t,\textbf{x})) = \sigma(t,\textbf{x}) \sigma(t,\textbf{x})^\top.
\end{equation}

The Pardoux-Peng theory \cite{pardoux-peng} gives a probabilistic representation, a nonlinear Feynman-Kac formula, for the terminal value problem solution $u(t,\textbf{x}),$ $t \in [0,T]$,
$\textbf{x}\in\mathbb{R}^{d}$ for the  $d$-dimensional quasi-linear parabolic PDE \eqref{eq-paraPDE1}.
Specifically, 
let $(X, Y, Z): [0,T]\times \Omega \to \R^d \times \R \times \R^d$ be an $\mathbb{F}$-adapted solution of the following FBSDEs:
\begin{align}%
dX_{t}  &  =\mu(t,X_{t},Y_{t})dt+\sigma(t,X_{t},Y_{t}%
)dB_{t}, \qquad
X_{0}    = \textbf{x}_0,
\label{fsde1}\\%
dY_{t}  &  =f(t,X_{t},Y_{t})dt+Z_{t}^{\top}\sigma(t,X_{t},Y_{t}%
)dB_{t}, \quad
Y_{T}    =g(X_{T}).
\label{bsde1}%
\end{align}
\index{Pardoux--Peng theory of nonlinear Feynman--Kac formula}
It was shown in \cite{pardoux-peng} that if $u \in C^{1,2}$ solves \eqref{eq-paraPDE1} and \eqref{termData}, then $(X,Y,Z)$ can be represented in terms of $u$ along the forward diffusion $X_t$, i.e.,
\begin{equation}
Y_{t}=u(t,X_{t}),\quad Z_{t}=\nabla u(t,X_{t}), \quad t \in [0, T),   \label{SDE-PDE}%
\end{equation} where the first identity in fact justifies the BSDE \eqref{bsde1} by using the It\^o formula \eqref{ItoForm2} on $u(t,X_{t})$  and the PDE \eqref{eq-paraPDE1}.

\subsection{Deep BSDE method}
The Deep BSDE method \cite{han2016deep,Han2018DeepBSDE} is based on the Pardoux--Peng theory in \Cref{sect_pardoux_peng} above.
It implements the discretization of the FBSDE \eqref{fsde1}--\eqref{bsde1} with a series of DNNs
\begin{equation}
f^{(n)}_{\theta^{n}}(\cdot) \approx \sigma(t_n, \cdot, u(t_n, \cdot))^\intercal \nabla u(t_n, \cdot)
\end{equation}
on time steps $t_n$, respectively.
These networks, together with $Y_0 = u(0, \mathbf{x}_0)$ as the initial value of the process $Y_t$, are trained using a loss function that matches the terminal value \eqref{termData}.
For a given time interval $[0,T]$, we define a partition%
\begin{equation}
    0=t_{1}<t_{2}<\cdots<t_{i}<t_{i+1}<\cdots<t_{N}=T,
    \label{T_partition}
\end{equation}
and, on each interval $[t_{n},t_{n+1}]$, use the Euler--Maruyama scheme for \eqref{fsde1} \eqref{bsde1}
to calculate $X_{n+1} \in \mathbb{R}^d$ and $Y_{n+1} \in \mathbb{R}$ with realizations of $B_n$,  i.e.,
\begin{align}
X_{n+1} &= X_{n}+\mu(t_{n},X_{n},Y_{n})\Delta t_{n}+\sigma(t_{n},X_{n},Y_{n})\Delta B_{n}, \label{euler1} \\
{Y_{n+1}} &= Y_{n}+f(t_{n},X_{n},Y_{n})\Delta t_{n}+ f^{(n)}_{\theta^{n}}(X_{n})^\intercal \Delta B_{n}, \label{euler2}
\end{align}
where $X_0=\textbf{x}_0 \in \mathbb{R}^d$, $\Delta t_{n}=t_{n+1}-t_{n}$, and $\Delta B_{n}=B_{n+1}-B_{n}$.
The output $Y_{N}$ is expected to be an approximation of $u(T,X_{N})$.
Thus, the loss function using the terminal condition is defined by a Monte Carlo approximation of
\begin{equation}
\mathcal{L}_{\rm BSDE}(\theta^{0}, \cdots, \theta^{N-1}, Y_0)= \mathbb{E}\left[ \left\|  Y_{N}-g(X_{N})\right\| ^{2} \right],
\end{equation}
whose minimizer $Y_0^{\ast}$ produces an approximation to $u(0,\textbf{x}_0)$ of \eqref{eq-paraPDE1} \eqref{termData}.

In this paper, to reduce the total number of DNNs and to solve the terminal value problem in a broader domain at $t=0$, the implementation of the Deep BSDE method utilizes two DNNs, one approximating the PDE solution $u(0, \cdot)$ in a region where $X_0$ samples from, the other one approximating $\nabla u(\cdot, \cdot)$ with continuous time input on $[0, T]$.


\subsection{Derivative-free martingale neural network - DeepMartNet} Here we construct the loss function for the neural network approximation $u_{\theta}(t,\textbf{x})$ by using
 Varadhan's  martingale problem formulation with diffusion processes \cite{Varadhan71}  for the solution of PDE (\ref{eq-paraPDE1}) \cite{Cai2023DeepMart,Cai2025SOCMart,Cai2025Mart}.
\index{DNN methods for PDEs!martingale based learning method}
\index{numerical methods!martingale based learning method}
For quasi-linear PDEs, we employ two stochastic processes for the training of the DNN $u_{\theta}(t,\textbf{x})$.

\smallskip
\noindent {\bf 1. (Pilot process)} Let $\hat{X}: [0, T] \times \Omega \to \R^d$ be a pilot process governed by
\begin{equation}\label{eq_SDE}
\begin{aligned}
    \hat{X}_t = \hat{X}_0 \;& + \int_0^{t} \mu\vbr{\big}{s, \hat{X}_s, u_{\bar \theta}(s, \hat{X}_s)} \di s + \int_{0}^t \sigma\vbr{\big}{s, \hat{X}_s, u_{\bar \theta}(s,\hat{X}_s)} \di B_s.
\end{aligned}
\end{equation}
The pilot process $\hat{X}_t$ is used to explore the state space $\R^d$; its sample paths do not need to be updated frequently during DNN training by freezing NN's parameters in $u_{\bar \theta}$ for many epochs between updates.


\noindent {\bf 2. (System process)} Let $h>0$, and 
for any $s \in [0, T-h]$, define a system process $\{X_t^s\}_{t\in[s, s+h]}$ starting from $\hat{X}_s$ (the state of the pilot process at time $s$) and generated by the operator $\mathcal{L}$, namely, for $t \in [s, s + h]$,
\begin{equation}\label{eq_branch}
    X_t^s = \hat{X}_s + \int_s^t \mu\br{r, X_r^s, u_{\theta}(r, X_r^s)} \di r + \int_s^t \sigma\br{r, X_r^s, u_{\theta}(r, X_r^s)} \di B_r,
\end{equation}
where the superscript $^s$ indicates that $X_t^s$ starts at the time $s$.

By It\^o's formula \eqref{ItoForm2}, $u_{\theta}(t,\vx)$ will exactly satisfies the PDE~\eqref{eq-paraPDE1} for $(t,\vx)$ in the region explored by $\hat{X}_t$ with positive probability, 
provided the following martingale condition holds \cite{Varadhan71,Cai2023DeepMart}:
\begin{equation}\label{eq_martPDE}
    \E{\M_{t+h}^t \big\vert \hat{X}_t} = 0, \quad 0 \leq t \leq T-h,
\end{equation}
where the increment process $\M_t^s$ is defined by
\begin{equation} \label{eq_defMt} 
    \M_t^s := u_{\theta}(t, X_t^s) - u_{\theta}(s, X_s^s) + \int_s^t f(r, X_r^s, u_{\theta}(r, X_r^s)) \di r, \;\; t \in [s, s + h]. 
\end{equation} 

As indicated by \eqref{eq_martPDE}, the pair $(X_{t+h}^t,\M_{t+h}^t)$ captures local information of the PDE~\eqref{eq-paraPDE1} at $(t,\x) = (t,\hat{X}_t)$. Therefore, their sample paths should be re-calculated during training to reflect the updates in $u_\theta$.
However, the calculation is only required over the short interval $[t,t+h]$, and these local simulations are parallelizable for different starting times $t$, a key feature of this method.

\bigskip
\noindent{\bf DeepMartNet  in Galerkin weak formulation}
To enforce the martingale condition in \eqref{eq_martPDE} without computing the conditional expectation, we introduce a test function $\rho(t, \hat{X}_t)$; it can be shown \cite{Cai2025SOCMart,Cai2025Mart} that \eqref{eq_martPDE}  is equivalent to a Galerkin weak form
\begin{equation}
    \int_0^{T-h} \E{\rho(t, \hat{X}_t) \M(t, \hat{X}_t, \xi; u_{\theta})} \di t =0 + \mathcal{O}(h^2). 
\end{equation}
where the expectation is taken with respect to the PDF of $(\hat{X}_t, \xi)$, and 
\begin{equation}\label{eq_defDeltaM}
    \begin{aligned}
        \M\br{t, \textbf{x}, \xi; u} :=\;& u\br{t+h, \textbf{x} + \mu(t,\textbf{ x}, u(t,\textbf{ x})) h + \sigma(t, \textbf{x}, u(t, \textbf{x})) \sqrt{h} \xi} \\
        &-u(t, \textbf{x}) + h f(t, \textbf{x},u(t,\textbf{ x}))
    \end{aligned}
\end{equation}
for $(t, \textbf{x}, \xi) \in [0, T] \times \R^d \times \R^d$ and $\xi \sim \mathrm{N}(0, I_d)$. 
The remainder term $\mathcal{O}(h^2)$ comes from the discretization error of the integrals in \eqref{eq_branch} and \eqref{eq_defMt}.

Then, a completely derivative-free DeepMartNet is to solve the following mini-max problem for a NN approximation $u_{\theta}(t, \x)$ to the PDE solution $u(t, \x)$,
\begin{equation}\label{eq_wPMF}
    \min_{\theta}  \; \max_{\eta} \abs{G(u_{\theta}, \rho_{\eta})}^2, \;\; G(u_{\theta}, \rho_{\eta}) := \int_0^{T-h} \E{\rho_{\eta}(t, \hat{X}_t) \M(t, \hat{X}_t, \xi; u_{\theta})} \di t,
\end{equation}
where $\rho_{\eta}$ is the test function DNN.
In practice, the minimax problem \eqref{eq_wPMF} will  be solved via adversarial training \cite{Zang2020Weak}; see \cite{Cai2025SOCMart,Cai2025RDO} for details of implementation.

\smallskip
\noindent{\bf Remark} (Increment $\M_{s+h}^s $ and PDE residual)
By the It\^o formula \eqref{ItoForm2}, we can find a connection between the process $\{\M_t^s\}_{t\in[s, s+h]}$ and PDE residual:
\begin{equation}\label{eq_EMeqhR}
    \E{\M_{s+h}^s \vert \hat{X}_s} = h R(s, \hat{X}_s; u_{\theta}) + \mathcal{O}(h^2), \quad 0 \leq s \leq T-h,
\end{equation}
where $R(t, \textbf{x}; u_{\theta}) := (\partial_t + \mathcal{L}) u_{\theta} (t, \textbf{x}) -f(t, \textbf{x}, u_{\theta}(t, \textbf{x}))$ denotes the residual error of \eqref{eq-paraPDE1} for the neural network approximation $u_{\theta}$. Therefore, the martingale neural network can be viewed as an equivalent to a Galerkin implementation, weighted with the PDF of the pilot paths, of a PINN method. 

For implementation details of DeepMartNet,  refer to Appendix A.2

\subsection{Random difference and shotgun neural networks} 
Equation \eqref{eq_defDeltaM} can be viewed as a random difference method (RDM) for the PDE operator as seen from \eqref{eq_EMeqhR}. A general approach of deriving such RDMs can be done using simple Taylor expansions \cite{Xu2025Shotgun,Cai2025RDO}, which were in fact also used to derive the It\^o formula  \eqref{ItoForm2} \cite{weinanbook}.

Following \cite{Cai2025RDO}, 
we let $\xi = \br{\xi_1, \xi_2, \cdots, \xi_d}^{\top}$ be a $d$-dimensional random vector. 
By the Taylor expansion of a function $f(\z)$ around $\z = 0$, we obtain
\begin{equation}\label{eq_taylor}
    \E{f(\sqrt{h} \xi)} = f(0) + \sum_{k=1}^3 \frac{h^{k/2}}{k!}\E{\vbr{\big}{\xi^{\top} \partial_{\z}}^k f(0)} + \frac{h^2}{4!}\E{\vbr{\big}{\xi^{\top} \partial_{\z}}^4 f(c\sqrt{h} \xi)},
\end{equation}
where $\xi^{\top} \partial_{\z} := \sum_{i=1}^d \xi_i \partial_{z_i}$,
and $c \in [0, 1]$ depending on $\sqrt{h} \xi$.
Assume 
\begin{equation}\label{eq_condi_xi}
    \E{\xi_i} = 0, \quad \E{\xi_{ij}} = \delta_{ij}, \quad \E{\xi_i \xi_j \xi_k} = 0, \quad \E{\abs{\xi_i \xi_j \xi_k \xi_l}} < \infty
\end{equation}
for $i, j, k, l = 1, 2, \cdots, d$, 
where $\delta_{ij}$ denotes the Kronecker delta. 
Under these conditions, the first- and third-order derivative terms in \eqref{eq_taylor} vanish, yielding
\vspace*{-0.6cm}
\begin{equation}\label{eq_rdo_expa}
    \E{f(\sqrt{h} \xi)} = f(0) + \frac{h}{2} \sum_{i=1}^d \partial_{z_i}^2 f(0) + \mathcal{O}\br{h^2}.
\end{equation}

Next, we consider the operator $ \mathcal{D} = \partial_t + \mathcal{L}$
appearing in the parabolic PDE \eqref{eq-paraPDE1}.
For a given point $(t, \textbf{x}) \in [0, T] \times \R^d$, consider function $U(s,\textbf{z})$ given by
\begin{equation}\label{eq_defVsz}
    U(s, \textbf{z}) := u(t + s, \textbf{x} + \mu s + \sigma \textbf{z}), \quad (s, \textbf{z}) \in [T - t] \times \R^d,
\end{equation}
for $u \in C^{2,4}([0, T] \times \mathbb{R}^d)$.
Applying \eqref{eq_rdo_expa} with $f$ replaced by $U(h, \cdot)$, we have
\begin{equation}\label{eq_expandV}
    \E{U(h, \sqrt{h}\xi)} = U(h, 0) + \frac{h}{2} \sum_{i=1}^d \frac{\partial^2 U}{\partial z_i^2} (h, 0) + \mathcal{O}\br{h^2}. 
\end{equation}
By the definition in \eqref{eq_defVsz}, the Taylor expansion and the chain rule imply
\begin{align}
   U(h, 0) &= U(0, 0) + \br{\partial_t u(t, \x) + \mu^{\top} \nabla u(t, \x)} h + \mathcal{O}(h^2), \label{eq1_expandV}\\
    \sum_{i=1}^d \frac{\partial^2 U}{\partial z_i^2} (h, 0) &= \sum_{i=1}^d \frac{\partial^2 U}{\partial z_i^2}(0, 0) + \mathcal{O}(h) = \mathrm{Tr}\rbr{\sigma \sigma^{\top} \nabla\nabla^{\top} u(t, \x)} + \mathcal{O}(h). \label{eq2_expandV}
\end{align}
Inserting \eqref{eq1_expandV} and \eqref{eq2_expandV} into \eqref{eq_expandV}, we obtain 
\begin{equation}\label{eq_expandV2}
    \E{U(h, \sqrt{h}\xi)} = U(0, 0) + h \mathcal{D} u(t, \textbf{x}) + \mathcal{O}(h^2). 
\end{equation}
where $\mathcal{D}=\partial_t + \mathcal{L}$ as in the PDE \cref{eq-paraPDE1}. Thus, by setting
\begin{equation}
    D_h(t, \textbf{x}, \xi; u) = \left(u(t + h, \textbf{x} + \mu h + \sqrt{h} \sigma \xi) - u(t, \textbf{x})\right)/h,
\end{equation}
we have a random difference approximation by expectation
\begin{equation}\label{eq_rdo}
\mathcal{D}_h u(t, \textbf{x}) := \E{D_h(t, \textbf{x}, \xi; u)} = \mathcal{D}u(t, \textbf{x}) + \mathcal{O}(h).
\end{equation}
A rigorous derivation of \eqref{eq_rdo}, along with an explicit upper bound for the remainder term $\mathcal{O}(h)$, is provided in \cite{Cai2025RDO}.
The random difference approximation  has also been obtained for quasi-linear parabolic PDEs  in \cite{Xu2025Shotgun}, by single-step discretization of FBSDE nonlinear Feynman-Kac formula \cite{Raissi2023FBSNNs}, and in the early work of DGM \cite{sirignano2018dgm}.
A feature of \cref{eq_rdo} is that it approximate 2nd order derivatives using only first-order differences, avoiding computing derivatives.

Approximating $\mathcal{D}$ with $\mathcal{D}_h$ in \eqref{eq_rdo}, the PDE \eqref{eq-paraPDE1} can be solved by training a neural network $u_\theta$ satisfying
\begin{equation}\label{eq_ER0}
    \E{R(t, \textbf{x}, \xi; u_{\theta})} = 0, \quad (t, \textbf{x}) \in [0, T-h] \times \R^d,
\end{equation}
where $R(t,\textbf{ x}, \xi; u_{\theta})$ is the residual function defined by
\begin{align}
    R(t, \textbf{x}, \xi; u_{\theta}) &= D_h(t, \textbf{x}, \xi; u_{\theta}) - f(t, \textbf{x}, u_{\theta}(t, \textbf{x})). \label{eq_defR}
\end{align}
In \cite{Cai2025RDO}, a Galerkin method is adopted to reduce the residual in \eqref{eq_ER0}, resulting in the Deep random difference method (DRDM),  an equivalent to DeepMartNet. 

The shotgun method \cite{Xu2025Shotgun} differs in its implementation details.
First, an averaged alternative of \cref{eq_rdo} using antithetic variates of the random vector $\xi$
\begin{equation}\label{eq_rdo_av}
\bar{\mathcal{D}}_h u(t,\textbf{ x}) :=\frac{1}{2} \E{D_h(t,\textbf{x}, \xi; u) + D_h(t, \textbf{x}, -\xi; u)}
\end{equation}
is applied to achieve a first order accuracy in $h$.
The shotgun method enjoys a simple loss function formulation as in the standard PINN method.


\section{Numerical results}

In this section, we will provide the numerical results of the neural network methods discussed in Sections 3 and 4 for two types of parabolic equations \eqref{eq-paraPDE1} in $\mathbb{R}^d$, with dimensionality $d=100,  1000$.
All codes and input files for producing the numerical results can be found at \url{https://github.com/caiweismu/DNN_HighDimPDEs}.


\subsection{Hamilton-Jacobi-Bellman (HJB) equation}  We consider a HJB equation with a linear-quadratic-Gaussian (LQG) control as follows:
\begin{equation}\label{eq_hjb}
\partial_t u( t,\textbf{x}) + \frac{1}{2}\Delta u(t,\textbf{x}) - \Vert \nabla_{x} u(t,\textbf{x}) \Vert^2 = 0,
\end{equation} 
where $(t,\textbf{ x}) \in [0, T) \times \R^d$, 
with the terminal condition $u(T, \textbf{x}) = g(\textbf{x}),\quad T = 1$ and
 a non-convex Rosenbrock function for $c_{1,i},c_{2,i} \sim \text{Unif}[0.5, 1.5]$
\begin{equation}\label{eq_rosenbrock}
g(\textbf{x}) = \log\left(\frac{1+\sum_{i=1}^{d-1}\left[c_{1,i}( x_{i} - x_{i+1})^2 + c_{2,i}x_{i+1}^2\right]}{2}\right).
\end{equation}

By the Cole-Hopf transformation, the value function $u(t, \textbf{x})$ is given by 
\begin{equation}
\label{hjb_sol}
    u(t,\textbf{x}) = -\frac{1}{2}\log\left(\int_{\mathbb{R}^d}(2\pi)^{-d/2}\exp(-\Vert \boldsymbol{y} \Vert^2  /2)\exp(- 2 g(\textbf{x} - \sqrt{(1-t)}\boldsymbol{y}))d\boldsymbol{y}\right).
\end{equation}
For reference solutions, we compute $u(t, \textbf{x})$ by approximating the expectation using the Monte Carlo method, where $10^6$ samples of unit Gaussian are applied to \eqref{hjb_sol}.
The neural network solution will be tested on points sampled from $\textbf{x} \sim \mathcal{N}(0, I_d)$, $t=0$. Table \ref{hjbtable} lists the performance of the tested neural networks whose network structures are given in Appendix A.4.1.

\begin{table}
    \centering
    \caption{ Numerical results for HJB equation \eqref{eq_hjb}  on points sampled from $\textbf{x} \sim \mathcal{N}(0, \boldsymbol{I}), t=0$. $\text{RE}_2$ - relative $L^2$ error, MB = MB in Memory, Time = GPU time in seconds (s) or minutes  (m).}
    {\small
    \begin{tabular}{lcccc}
      \toprule
         &PINN-SDGD  & Deep BSDE  & DeepMartNet  & Shotgun \\
         \midrule
     Network size    & $4 \times 1024$ & $(4 \times 1024)\times 2$ & $4 \times 1024$ & $4 \times 1024$ \\
     \midrule
      $\text{RE}_2$ (d=100)  & $3.12 \times 10^{-3}$ & $5.02 \times 10^{-3}$  & $1.35 \times 10^{-2}$ & $5.22 \times 10^{-3}$  \\
           \midrule
      $\text{RE}_2$ (d=1000)  & $6.08 \times 10^{-3}$ & $3.76 \times 10^{-3}$& $3.14 \times 10^{-3}$ & $2.59 \times 10^{-3}$\\
      \midrule
       MB (d=100)   & 1127 & 1792 & 2273  & 2654 \\
        \midrule
       MB (d=1000)   & 1844 & 1856 & 9633  & 8192\\
       \midrule
      Time (d=100)   &  185m (A100)& 197s (4090) & 216s (A100) & 288s (4090) \\
       \midrule
      Time(d=1000)    & 217m (A100) & 678s (4090)& 260s (A100) & 867s (4090)\\
      \bottomrule
    \end{tabular}
    }
        \label{hjbtable}
\end{table}

\subsection{Black-Scholes (BS) equation}
%
Consider the BS equation \cite{beck2021solving} 
\begin{equation} \label{eq_bse}
u_t(t,\textbf{x}) + \frac{1}{2}\sum_{i=1}^d|\sigma_ix_i|^2\frac{\partial^2 u(t,\textbf{x})}{\partial x_i^2} + \sum_{i=1}^d\mu_ix_i\frac{\partial u(t,\textbf{x})}{\partial x_i} = 0,
\end{equation}
whose corresponding SDE trajectory is a geometric Brownian motion (GBM). {\color{black}
We set the risk‐free interest rate $r=1/20$, the drift parameter $\mu=-1/20$, and
the volatility is set to be $\sigma_i = 0.1 + 0.4 i/d,  i=1, 2,\ldots,d.$
}
The terminal condition at $T=1$ encodes the contract agreement by 
\begin{equation}\label{eq_bse_terminal}
u(T,\textbf{x}) = \phi_{100}(\textbf{x}), \qquad  \phi_{a}(\textbf{x}) = \exp(-rT)\max\left( {\max_i x_i - a,0}\right).
\end{equation}
The exact solution can be simulated by a Monte Carlo method:
\begin{equation}\label{eq_exact_bse_corr}
u(t,\textbf{x}) = \mathbb{E}\left[\phi_{100}\left(y_1, y_2, \cdots, y_d\right)\right], y_i = x_i \exp\left(\sigma_i  W_{T-t, i} + (\mu_i - \sigma_i ^2/2)(T-t)\right). 
\end{equation}
We examine the solution at $t=0$ within $x \in [90, 110]^d$ and at $t=0$ for $d =100, 1000$ at 1000 test points uniformly from the cube  \cite{beck2021solving}. The results are summarized in Table 2 with the network structures given in Appendix A.4.2.
For SDE based methods, a scaled version of \eqref{eq_bse} is used with the following scaling $\textbf{z} := \textbf{x} / K$, $v(t,\textbf{z}) := u(K\textbf{ z}, t) / K$, $K=100$.
The scaled BS equation takes the exact same form of \eqref{eq_bse} with $\textbf{x}$ replaced by $\textbf{z}$
and $u$ by $v$ and the exact solution is given by $v(t,\textbf{z}) = \mathbb{E}\left[\phi_{1}\left(y_1, y_2, \cdots, y_d\right)\right].$

\begin{table}
    \centering
    \caption{ Numerical results for BS equation \eqref{eq_bse}  on uniform points in $[90, 110]^d$. $\text{RE}_2$ - relative $L^2$ error, MB = MB in Memory, Time = GPU time-seconds (s). }
    {\small
    \begin{tabular}{lcccc}
      \toprule
         &PINN  & Deep BSDE  & DeepMartNet  & Shotgun \\
         \midrule
     Network size    & $4 \times 128$ & $(4 \times 1024)\times 2$ & $4 \times 1024$ & $4 \times 1024$ \\
     \midrule
      $\text{RE}_2$ (d=100)  & $1.42 \times 10^{-2}$ & $1.35 \times 10^{-2}$  & $1.59 \times 10^{-2}$ & $1.95 \times 10^{-2}$  \\
           \midrule
      $\text{RE}_2$ (d=1000)  & $4.69 \times 10^{-3}$ & $6.17 \times 10^{-3}$& $8.27 \times 10^{-3}$ & $1.48 \times 10^{-2}$\\
      \midrule
    MB  (d=100)   & 510 & 1785 & 1334  & 2027 \\
        \midrule
    MB (d=1000)   & 3174 & 1896& 8310  & 7431\\
       \midrule
      Time (d=100)   & 60s (T40) & 185s (4090) & 93s (A100) & 170s (4090) \\
       \midrule
      Time (d=1000)    & 276s(T40) & 676s(4090)& 110s (A100) & 642s(4090)\\
      \bottomrule
    \end{tabular}
    }
    \label{bstable}
\end{table}

\begin{table}[!htbp]
\centering
\caption{Comparison of PDE and SDE-based and hybrid DNNs for solving high-dimensional quasi-linear parabolic equation \eqref{eq-paraPDE1} and BC \eqref{BC}.}
{\small
\begin{tabular}{|p{1.4cm}|p{2.8cm}|p{2.8cm}|p{3.2cm}|}\hline
DNN sol. $u_\theta(t, \textbf{x})$ & {\bf Strategy A: PDE based DNNs} & {\bf Strategy B: SDE based DNNs}& {\bf Strategy C: Random differences}\\
\hline
Math. behind the loss& PDE residuals at $\textbf{t}=\{t_j\}, \textbf{X}=\{\textbf{x}_i\}$,  $R(u_\theta,\textbf{t},\textbf{X})$ and BC& Probabilistic form for PDE solution and BC& It\^o formula \eqref{ItoForm2} applied on diffusion paths and BC\\
\hline
Sampling strategy& Monte Carlo sampling of domain by a given distribution or SDEs &  Random walk by stochastic systems with  the PDE operator generator &  Random walk by stochastic systems with the PDE operator generator\\\hline
$\nabla$, $\nabla \nabla^T$&         SDGD, HTE, RS &   It\^o formula \eqref{ItoForm2}&  It\^o formula \eqref{ItoForm2}\\\hline Training aims& Reducing PDE residuals of DNNs at sampled locations in strong form
in $l_2$ or $l_{\infty}$ 
& Enforcing the nonlinear Feynman-Kac formula \eqref{SDE-PDE} for PDE solution and FBSDEs&  Reducing the PDE residual of neural network in $\pi(X_t)$-weighted  Galerkin weak form\\
&  \vspace*{0pt}  $\min_\theta||R(u_\theta,\textbf{t},\textbf{X})||$& $u_\theta(t,X_t)=Y_t$, $\nabla u_\theta(t, X_t)=Z_t$&$\min_\theta \max_\eta$ $(\rho_\eta,R(u_\theta, t, X_t))_{\pi(X_t)}$\\ \hline
Pros& Parallel, any PDEs, error estimate for elliptic PDEs, SDGD-PINN/RS-PINN for derivatives& Hessian matrix-free, system path-importance sampling, parallel between SDE paths&Total derivative-free, path-importance sampling, parallel within each path,  first order accuracy in $\Delta t$, reduced loss variance.\\\hline
Cons& Simple PINNs require derivatives, lack sampling distribution& need paths, not parallel within a path, only parabolic PDEs&need paths, test network $\rho_\eta$, min-max saddle point, only parabolic PDEs\\\hline
      Rep. works&PINN \cite{raissi2019pinn}, DeepRitz \cite{Weinan2017TheDR,YingSemigroup22}, DGM \cite{sirignano2018dgm}, SDGD \cite{hu2023tackling}, HTE \cite{hu2023hutchinson}, RS-PINN \cite{he2023learning}&  DeepBSDE \cite{Han2018DeepBSDE}, FBSNN \cite{Raissi2023FBSNNs,Zhang2022FBNN}, Actor-Critic \cite{Zhou2021Actor}, Diffusion \cite{Nusken2023JML,Ruthotto24}, Feynman-Kac \cite{HanNica2020} \cite{YingSemigroup22}, Picard Iter\cite{han2024picard}, Diffusion MC \cite{DFMC20}&  Weak form using test function: DeepMartNet \cite{Cai2025SOCMart,Cai2025Mart}, DRDM \cite{Cai2025RDO} and, strong PINN-form: Shotgun \cite{Xu2025Shotgun}, DGM \cite{sirignano2018dgm}, RS-PINN \cite{he2023learning}\\ \hline
\end{tabular}
}
\label{nn-overview}
\end{table}

\section{Conclusion and outlook}
To find efficient and accurate solutions to very high-dimensional PDEs, the following key issues have to be addressed: (1) strategy of sampling the solution domain; (2) computation and treatment of the solution derivatives in high-dimensional spaces; (3) bias and variance reduction of the mini-batch implementation of the loss function; (4) parallel implementation of the training algorithms such as the computation of loss functions, etc; (5) the demand for excessive memory. Table \ref{nn-overview} summarizes the attributes of a  list of selected neural network methods developed in recent years.

\subsection{Strength and weakness of discussed methods} Among the four methods presented in this paper (PINN, Deep BSDE, DeepMartNet, and Shotgun/DRDM),  we find that  PINN has the simplest way of (parallel) implementation, and, with the SDGD and RS-variants, it provides accuracy and efficiency in the sampled solution domain of time-space, and  can handle the high-dimensional PDEs tested in this work. However, it does experience some challenges to provide accurate approximation to solution with large gradients, and hence research on more effective sampling to capture the large changes of solution is an important future research area. One possibility is to leverage the SDE sampling strategy to improve the sampling for the PINN type methods.

DeepBSDE has proven to be an effective method to find the solution of high dimensional parabolic PDEs  at $t=0$, especially, when the original DeepBSDEs are augmented with two neural networks to represent the solution at $t=0$ and the gradient for all time, respectively. The method is parallel in the sampling of different paths of the SDEs, however, not parallelizable within each path. The martingale neural network DeepMartnet has been shown to be robust and effective, and  enjoys parallelism even within each SDE path, reduced variance in the loss function calculation, and can resolve high gradients for very high-dimensional PDEs and stochastic optimal controls problems in up to 10,000 dimensions \cite{Cai2025SOCMart,Cai2025Mart}. 
The shotgun method presents a PINN-type random finite difference method, similar to RS-PINN while the DRDM is its weak form implementation, an equivalent to the DeepMartNet. The random difference methods provide approximation to PDE solutions in regions explored by the SDE paths, in particular, with a better accuracy near $t=0$.
\subsection{Open problems} Further research will be needed to address the following outstanding issues for the methods discussed in this paper. 

 \smallskip
 {\textbf{1. PDE-residual based methods}.}  \textbf{$\bullet$Applicability} \textcolor{black}{Compared to other methods, PINNs stand out for their generality, being applicable to all types of PDEs, particularly those with complex boundary conditions and irregular domains.}
\textbf{$\bullet$Algorithm bottlenecks} \textcolor{black}{The main bottleneck of PINNs lies in the sampling of residual/collocation points. Although adaptive and effective sampling methods have been well-studied for low-dimensional PDEs, typically on bounded domains, this remains under-explored for high-dimensional, unbounded domains. Sampling along SDE trajectories works well in some cases but may fail for solutions with complex structures and high gradients in very high dimensions, especially for time-dependent problems requiring long-time integration. The key challenge is the lack of a theoretical guarantee that minimizing the residual plus initial-loss along SDE trajectories will ensure effective approximation of the true solution. Nonetheless, similar theories for bounded domains have been extensively studied \cite{hu2021extended,Ryck2021ErrorAF,mishra2020estimates}. Thus, how to design solution-adapted sampling strategies will be a research topic of urgent need.}  \textbf{$\bullet$Computational bottlenecks} \textcolor{black}{An open problem is how to avoid the high cost of gradient and Hessian matrix calculation for variable coefficients and nonlinear differential operators by extending the SDGD and RS-PINN. More effective variance reduction strategies in the mini-batch calculation of the PINN loss function for high-dimensional PDEs should also be explored.}

\smallskip
 {\textbf{2. SDE based methods}.} \textbf{$\bullet$Applicability} As the SDE based NN requires probabilistic representation for the PDEs' solution, this type of methods are limited to parabolic and elliptic problems on either bounded domain or whole space. An open problem is the treatment of boundary conditions in bounded domains, especially the Neumann and Robin boundary conditions. \textbf{$\bullet$Algorithm bottlenecks} Most of the SDE based methods can only provide solution to parabolic PDEs  in a low-dimensional manifold at $t=0$ or sometimes in the limited region explored by the sampling stochastic processes as shown in the case FBSNN \cite{Raissi2023FBSNNs}. How to extend the capability of these type of methods to the whole solution time-space domain will be very challenging. \textbf{$\bullet$Computational bottlenecks} Parallel processing of data within a sampling path in computing the loss function  and ways of variance reduction in the mini-batch computation of the loss function should also be studied.

\smallskip
 {\textbf{3. Hybrid methods - random differences}} 
 \textbf{$\bullet$Applicability} The DeepMartNet and random difference methods are limited to parabolic and elliptic PDEs due to SDE approaches involved.\textbf{$\bullet$Algorithm bottlenecks} Hybrid methods, sharing features of two approaches above, only find good solution of parabolic PDEs in a local region near $t=0$. Moreover, as the adversarial training is required for the Galerkin approach, robustness of solving the min-max saddle point needs to be improved in the stochastic gradient methods.
  \textbf{$\bullet$Computational bottlenecks} The cost comes from sampling SDE paths, especially the required updating of paths for quasi-linear PDEs.

\subsection{Outlook}
The future of neural network research for high-dimensional PDEs holds great promise by the fast developments of many existing methods within such a short time period, and their capability in obtaining non-trivial solutions for high-dimensional PDEs. The revolutionary technique of deep neural networks has engendered great hope in the computational community for cracking the decades-old hard problem of curse of dimensionality encountered in the solution of high-dimensional PDEs. The potential impact of new breakthroughs in neural network algorithms is far-reaching for the resolution of many outstanding challenges in scientific, engineering, and economic computing such as the Schr\"odinger equation for many-electron quantum systems, Hamilton-Jacobi-Bellman equations for stochastic optimal controls, and Fokker-Planck equations for complex biological systems and generative AIs, etc.

\bibliographystyle{siamplain}
\bibliography{ref}

\newpage 
\section*{Appendix: Algorithmic Details and Proof Sketches}

\subsection*{A.1 Notations} 

\begin{longtable}{ll}
\toprule
Symbol & Meaning \\
$u_\theta(t, \textbf{x})$ & Neural network approximation of $u(t, \textbf{x})$\\
$\nabla\nabla^\top u$ & Hessian matrix for function $u$\\
$R(u)$ & PDE residual for function $u$\\
$D = \partial_t + L$ & Parabolic operator \\
$\mu(t, \textbf{x}, u)$ & Drift coefficient \\
$\sigma(t, \textbf{x}, u)$ & Diffusion matrix \\
$A = \sigma\sigma^\top$ & Diffusion covariance \\
$\mathrm{Tr}(\cdot)$ & Trace operator \\
$B_t$ & Brownian motion \\
$\mathcal{F}_t$ & Natural filtration generated by $B_t$ \\
$X_t$ & Forward SDE (state process) \\
$Y_t$ & Backward SDE (value process) \\
$Z_t$ & Backward SDE (gradient process) \\
$M_t^{s}$ & Martingale increment from time $s$ to $t$ \\
$\rho_\eta(t, \textbf{x})$ & Test function in adversarial training\\
\bottomrule
\end{longtable}

\subsection*{A.2 Implementation of DeepMartNet} For implementations of DeepMartNet, we introduce a uniform time partition on $[0, T]$, i.e., $t_n = nh$ for $n = 0, 1, \cdots, N$, with step size $h = T/N$.
Then, the loss function in \eqref{eq_wPMF} can be replaced by its mini-batch version, i.e.,
\begin{align}
    &\abs{G(u, \rho)}^2 \approx G^{\top}(u, \rho; A_1) G(u, \rho; A_2), \label{eq_approxG2}\\
    &G(u, \rho; A_i) := \frac{h}{\abs{A_i}} \sum_{(n, m) \in A} \rho(t_n, \hat{X}_{n}^m) \M(t_n, \hat{X}_{n}^m, \xi_n^m; u), \quad i = 1, 2,\label{eq_defGA}
\end{align} 
where $\abs{A_i}$ denotes the number of elements in the index set $A_i$; $\{\hat{X}_{n}^m\}_{n=0}^N$ are the pilot paths of $\hat{X}$ generated by applying the Euler-Maruyama scheme to \eqref{eq_SDE}, i.e., 
\begin{equation}\label{eq_iidXxi}
    \hat{X}_{n+1}^m = \hat{X}_n^m + \mu\br{t_n, \hat{X}_n^m, \hat{u}(t_n, \hat{X}_n^m)}\,h + \sigma\br{t_n, \hat{X}_n^m, \hat{u}(t_n, \hat{X}_n^m)}\,\sqrt{h}\,\xi_{n}^{m}
\end{equation}
for $n = 0, \cdots, N-1$ and $m = 1, \cdots, M$, where each $\xi_{n}^{m}$ is drawn i.i.d. from $\mathrm{N}(0, I_q)$.
The index sets $A_1$ and $A_2$ are constructed by 
\begin{equation}\label{eq_defA1A2}
    A_i = \mathrm{N}_i \times \mathrm{M}_i, \quad \mathrm{N}_i \subset \{0, 1, \cdots, N-1\}, \quad \mathrm{M}_i \subset \{1, 2, \cdots, M\}, \quad \mathrm{M}_1 \cap \mathrm{M}_2 = \emptyset
\end{equation}
for $i = 1, 2$. 
Here, inspired by \cite{Guo2022Monte,hu2023tackling}, the disjointness between $\mathrm{M}_1$ and $\mathrm{M}_2$ is used to ensure the mini-batch estimation in \eqref{eq_approxG2} is unbiased. 
The sets $\mathrm{N}_1$ and $\mathrm{N}_2$ can be chosen randomly, wthout replacements, from $\{0, 1, \cdots, N-1\}$, or simply set to $\{0, 1, \cdots, N-1\}$ if $N$ is not too large. 
By the estimation in \eqref{eq_approxG2}, the minimax problem in \eqref{eq_wPMF} can be solved by adversarial training methods as in \cite{Zang2020Weak}, resulting in a derivative-free DeepMartNet.

\subsection*{A.3 PINN algorithm details}\label{appendix:pinn_details}
\subsubsection*{A.3.1 SDGD detailed loss function}
Under the SDGD notation, instead of optimizing the expensive full PINN residual loss
\begin{equation*}
\mathcal{L}_{\text{PINN}}(\theta) = \frac{1}{2}\left(\sum_{i=1}^d \mathcal{L}_iu_\theta(\vx) - f(\vx)\right)^2,
\end{equation*}
we plug the estimated residual into the residual loss with a randomly sampled index set $I \subset \{1,2,\cdots,d\}$:
\begin{equation*}
\mathcal{L}_{\text{SDGD-biased}}(\theta) = \frac{1}{2}\left(\textcolor{black}{\frac{d}{|I|}\sum_{i \in I}} \mathcal{L}_iu_\theta(\vx) - f(\vx)\right)^2. 
\end{equation*}
Due to the nonlinearity of the mean square error loss, it is biased. We can sample the index set twice with $I, J \subset \{1,2,\cdots,d\}$ independently to debias, as detailed in Appendix B.
\begin{equation*}
\mathcal{L}_{\text{SDGD-unbiased}}(\theta) = \frac{1}{2}\left(\textcolor{black}{\frac{d}{|I|}\sum_{i \in I}} \mathcal{L}_iu_\theta(\vx) - f(\vx)\right)
\left(\textcolor{black}{\frac{d}{|J|}\sum_{j \in J}} \mathcal{L}_ju_\theta(\vx) - f(\vx)\right).
\end{equation*}
In practice, we find that the bias of the biased loss $\mathcal{L}_{\text{SDGD-biased}}(\theta)$ is negligibly small to produce as good convergence as the unbiased one, while the biased loss only samples once and saves time.

\subsubsection*{A.3.2 HTE detailed loss function}
Under the HTE notation, the full expensive PINN loss is
\begin{equation}
L_{\text{PINN}}(\theta) = \frac{1}{2}\left(\mathrm{Tr}(A\nabla \nabla^\intercal u_\theta) - f\right)^2,
\end{equation}
with the biased and unbiased versions:
\begin{equation}\label{eq:HTE}
L_{\text{HTE-biased}}(\theta;\{\vv_i\}_{i=1}^V) = \frac{1}{2}\left(\frac{1}{V}\sum_{i=1}^V\vv_i^\intercal\left(A\nabla \nabla^\intercal u_\theta\right)\vv_i -f\right)^2,
\end{equation}
\begin{align}
L_{\text{HTE-unbiased}}(\theta;\{\vv_i,\hat{\vv}_i\}_{i=1}^V) &= \frac{1}{2}\left(\frac{1}{V}\sum_{i=1}^V\vv_i^\intercal \left(A\nabla \nabla^\intercal u_\theta\right)\vv_i -f\right)\\
&\qquad\left(\frac{1}{V}\sum_{i=1}^V\hat{\vv}_i^\intercal \left(A\nabla \nabla^\intercal u_\theta\right)\hat{\vv}_i -f\right),
\end{align}
where $\{\vv_i,\hat{\vv}_i\}_{i=1}^V$ are $i.i.d.$ samples from $p(\vv)$.
\subsubsection*{A.3.3 Taylor-mode AD}
Taylor mode avoids the conventional stacked backward mode AD, tailored for high-order derivatives.
High-level speaking, its acceleration comes from considering the structure of high-order derivatives for nested and recursive neural net functions via the Faa di Bruno formula.
It can also be understood as a high-order extension (e.g., second-order HVP) of the efficient forward AD in the first-order JVP.

\subsubsection*{A.3.4 SDGD and HTE variance comparison}
SDGD and HTE estimate the trace of 2nd-order derivatives
$
\operatorname{Tr}(A) = \sum_{i=1}^d A_{ii}.
$
SDGD's estimator is
$
\operatorname{Tr}(A) \approx \frac{d}{|I|}\sum_{i\in I}A_{ii},
$
whose variance comes from the variance between diagonal elements across different dimensions.
Meanwhile, HTE's estimator is
$
\operatorname{Tr}(A) \approx\frac{1}{V}\sum_{k=1}^V \vv_k^\mathrm{T} A\vv_k,
$ 
where $V$ is the HTE batch size and each dimension of $\vv_k \in \mathbb{R}^d$ is an $i.i.d.$ sample from the Rademacher distribution.
Its variance is
$
\frac{1}{V} \sum_{i \neq j}A_{ij}^2.
$
HTE's variance only comes from the off-diagonal elements.

If the diagonal elements are similar, i.e., the PDE exact solution is symmetric, then SDGD has low variance.
If the off-diagonal elements of the Hessian are zero, i.e., different dimensions do not interact, then HTE is exact.
If the scales of off-diagonal elements are much larger than the diagonal ones, then HTE suffers from huge variance.

\subsubsection*{A.3.5 RS-PINN details}
RS-PINN's gradient and Hessian concerning the input $\vx$ are:
\begin{equation}
\nabla u_\theta(\vx) = \mathbb{E}_{\delta \sim \mathcal{N}(0, \sigma^2\boldsymbol{I})}\left[\frac{\delta}{\sigma^2}f_\theta(\vx+\delta)\right].
\label{rs-pinn-grad}
\end{equation}
\vspace*{-0,2cm}
\begin{equation}
\nabla \nabla^\intercal u_\theta(\vx) = \mathbb{E}_{\delta \sim \mathcal{N}(0, \sigma^2\boldsymbol{I})}\left[\frac{\delta\delta^{\mathrm{T}}-\sigma^2\boldsymbol{I}}{\sigma^4}f_\theta(\vx+\delta)\right].
\label{rs-pinn-hessian}
\end{equation}
They can be simulated using Monte Carlo sampling for the expectation estimator, which calculates derivatives without the need for expensive automatic differentiation. The PINN loss is then used to solve the PDE.

Compared with SDGD and HTE using the original network as the surrogate model, RS-PINN's model is smoothed and cannot fit nonsmooth PDE solutions, thus suffering from more approximation error. Further, the optimized Taylor-mode AD in SDGD and HTE ensures faster convergence than RS-PINN.

\subsection*{A.4.Network parameters for numerical results}
\subsubsection*{A.4.1 HJB equation results in Table 1}

{\color{black}The following summarize the network structure of each methods tested

\begin{itemize}
    \item  {\bf DeepMartNet:} The DNN $v_{\theta}$ consists of 4 hidden layers with width $W = 1024$. Training paths are generated via \eqref{eq_iidXxi} with initial states $X_0^{m} \sim \text{Unif}[90, 110]^d$. The number of time steps is $N=100$.
    \item {\bf PINN} The neural network consists of 4 hidden layers with width $W = 1024$. Training paths are generated via \eqref{eq_iidXxi} with initial states $X_0^{m} \sim \text{Unif}[90,110]^d$. We train PINN with $L^\infty$ loss achieved via adversarial training \cite{wang20222} and use SDGD \cite{hu2023tackling} for acceleration, with hyperparameters totally following Section 5.3 in SDGD \cite{hu2023tackling}.
    \item {\bf Deep BSDE:} The DNNs for desired initial values $u(0, \cdot)$ and for $\nabla u(\cdot, \cdot)$ both consist of 4 hidden layers, each with width $W=1024$.
    The number of time discretization steps is $N=100$.
    \item {\bf Shotgun:} The DNN for approximating PDE solution $u(\cdot, \cdot)$ consists of 4 hidden layers, with width $W=1024$.
    The number of discretization steps for SDE trajectories is $N+1 = 21$.
    The local time step size is $h=10^{-5}$.
    The local sample size for PDE residual estimator \cref{eq_rdo_av} is $M=8(32)$ for $d=100(1000)$.
    To facilitate training, terminal condition $u(T, x) = g(x)$ is embedded into the network with weight $w(t) = t/T$ for $d=100$ and $w(t) = \max\{0, (Nt/T - N + 1)^3\}$ for $d=1000$.
\end{itemize}
}

\subsubsection*{A.4.2 BS equation results in Table 2}
{\color{black}The following summarize the network structure of each methods tested

\begin{itemize}
    \item  {\bf DeepMartNet:} The DNN $v_{\theta}$ consists of 4 hidden layers with width $W = 1024$. Training paths are generated via \eqref{eq_iidXxi} with initial states $X_0^{m} \sim \text{Unif}[90, 110]^d$. The number of time steps is $N=100$.
    \item {\bf PINN} The neural network consists of 4 hidden layers with width $W = 1024$. Training paths are generated via \eqref{eq_iidXxi} with initial states $X_0^{m} \sim \text{Unif}[90,110]^d$. Vanilla PINN is used as the training algorithm with full residual computation.
    \item {\bf Deep BSDE:} The DNNs for desired initial values $u(0, \cdot)$ and for $\nabla u(\cdot, \cdot)$ both consist of 4 hidden layers, each with width $W=1024$.
    The number of time discretization steps is $N=100$.
    \item {\bf Shotgun:} The DNNs for approximating PDE solution $u(\cdot, \cdot)$ consists of 4 hidden layers, with width $W=1024$.
    The number of discretization steps for SDE trajectories is $N+1 = 21$.
    The local time step size is $h=10^{-5}$.
    The local sample size for PDE residual estimator \cref{eq_rdo_av} is $M=8(32)$ for $d=100(1000)$.
    To facilitate training, terminal condition $u(T, x) = g(x)$ is embedded into the network with weight $w(t) = t/T$ for $d=100$ and $w(t) = \max\{0, (Nt/T - N + 1)^3\}$ for $d=1000$.
\end{itemize}
}

\end{document}